\documentclass[fleqn
,12pt]{article}
\usepackage[affil-it]{authblk}
\usepackage[english]{babel}
\usepackage{blindtext}
\usepackage{enumitem}
\usepackage{indentfirst}

\providecommand{\keywords}[1]{\textbf{\textit{Keywords:}} #1}

\usepackage{graphicx}
\usepackage{amsmath}
\usepackage{amsthm}
\usepackage{amssymb}
\usepackage{amsfonts}
\usepackage[hidelinks]{hyperref}
\usepackage{hyperref}
\hypersetup{colorlinks, linkcolor=blue, citecolor=blue, urlcolor=blue}

\usepackage{booktabs}

\usepackage{longtable}
\usepackage{cite}
\usepackage{mathrsfs}

\usepackage[title,titletoc,toc]{appendix}

\newtheorem{case}{Case}

\theoremstyle{definition} 

\usepackage{chngcntr}
\counterwithin{case}{section}
\counterwithin{subcase}{case}
\counterwithin{exmp}{section}
\counterwithin{rem}{section}
\counterwithin{red}{section}
\counterwithin{prof}{section}

\usepackage[margin=1in,left=1in,right=1in]{geometry}
\usepackage{subcaption}
\usepackage{float}
\usepackage{placeins}
\allowdisplaybreaks
\usepackage[final]{changes}
\usepackage{appendix}
\usepackage{cleveref}
\usepackage{longtable}
\usepackage{rotating}
\usepackage{array}
\usepackage{booktabs}

\title{Exact self-similar solutions for nonlinear coupled heat–elastic wave systems with temperature-dependent material properties}

\author[1]{Manjit Singh\thanks{corresponding author: manjitcsir@gmail.com}}
\author[2]{Radhika}

\affil[1,2]{%
    Yadavindra Department of Sciences, Punjabi University Guru Kashi Campus, Talwandi Sabo--151302, Punjab, India.}
\begin{document}
\maketitle
\begin{abstract}
 In this work, we study a coupled nonlinear thermoelastic system in which the
heat capacity, thermal conductivity, and elastic modulus are all
temperature-dependent. The thermal conduction equation is formulated in a
generalized radial geometry governed by a parameter $\nu$, which encompasses
planar, cylindrical, and spherical configurations as special cases.
A complete Lie symmetry classification of the admissible constitutive functions
is carried out, yielding three principal classes: exponential laws, power laws,
and a unified class that contains both as limiting cases. For each class, the
admitted symmetry algebras are determined, optimal systems of one-dimensional
subalgebras are constructed, and the governing partial differential equations
are reduced to coupled ordinary differential equation systems via similarity
transformations.
Exact closed-form solutions are obtained for each class using a combination of
integrating factor methods, flux-variable substitutions, and power-law and
logarithmic ans\"{a}tze. For the exponential class, the thermal field exhibits
a logarithmic similarity profile while the displacement field decays
algebraically, with the geometry parameter absorbed into the leading-order
transformation. For the power-law class, the thermal field follows a self-similar
power-law profile and the displacement field decomposes into distinct elastic,
self-similar, and thermoelastically driven modes. Travelling-wave solutions are
also constructed, with elementary forms recovered in limiting cases. The
influence of geometry, particularly the structural distinction between
cylindrical and other configurations, is identified and physically interpreted.
All solutions are discussed in terms of the interplay between thermal diffusion,
elastic stiffness, and thermoelastic coupling.

\end{abstract}

\section{Introduction}
The mutual interaction between heat conduction and mechanical deformation is
a fundamental feature of many problems in continuum mechanics
\cite{biot1956thermoelasticity,vel2002exact,hetnarski2009thermal}. When a solid
body is subjected to a thermal load, the resulting temperature gradients
generate internal stresses that drive elastic displacements; conversely,
rapid deformation produces local heating through the thermoelastic coupling
term \cite{boley2012theory}. This two-way feedback is described by a coupled system of partial
differential equations (PDEs) linking the temperature field $u(z,t)$ to the
displacement field $w(z,t)$. Classical linear thermoelasticity, as
formulated by Biot \cite{biot1956thermoelasticity} and surveyed by Nowacki \cite{nowacki1970problems},
assumes constant material parameters and admits Fourier-series or integral
transform solutions. However, in many physically relevant situations, the 
material coefficients strongly depend on the temperature. High-temperature aerospace applications, geothermal processes in porous rocks, and thermal therapies in biological tissues frequently involve material properties whose thermal and mechanical characteristics exhibit significant temperature dependence. Consequently, the assumption of constant conductivity, heat capacity, and elastic modulus may no longer be adequate for accurate thermoelastic modelling \cite{jin1998thermal,rossmann2014review,nobrega2024review}. In such regimes,
linear theory breaks down and the governing equations become nonlinear and
coupled in an essential way.

Exact analytical solutions of nonlinear PDE systems are of particular value:
they provide benchmark data for numerical solvers, reveal the dominant
physical balance governing a given regime, and expose the role of constitutive
nonlinearity in a transparent form that parameter studies cannot easily
replicate \cite{qin2004blow,esquivel2023global,rabie2025investigation,ismail2025construction,ismail2025novel,rabie2025thorough}.

Self-similar and travelling-wave solutions, more generally group-invariant solutions, constitute an important class of exact solutions for coupled thermoelastic systems, since they reduce the governing partial differential equations to ordinary differential equations in suitable invariant variables. The existence of such reductions depends critically on the constitutive functions $C(u)$, $K(u)$, $E(u)$. For arbitrary material laws, nontrivial similarity reductions are generally not available. Consequently, the determination of all constitutive classes admitting additional symmetry generators becomes a group-classification problem, whose solution provides the foundation for the systematic construction of invariant solutions \cite{popovych2004new,ivanova2006group,vaneeva2007enhanced,nauryz2026lie}. This is the approach adopted in the present work.
The governing system considered here takes the form
\begin{equation}\label{thermo_elastic_model:1}
    \begin{aligned}
C(u)\,u_t&=\frac{1}{z^\nu}\left(K(u)\,z^\nu u_z\right)_z+\gamma\, w_t,
\\
\rho \,w_{tt}&=\left(E(u)\,w_z\right)_z+\beta\, u_z .
\end{aligned}
\end{equation}
Here $u=\,u(z,t)$ is the temperature or thermal field, $w=\,w(z,t)$ is the displacement field, $E(u)$ is the temperature dependent elastic modulus, $\beta\,u_{z}$ is the thermoelastic forcing caused by temperature gradients, the term $\left(E(u)\,w_z\right)_z$ implies that the local stiffness changes with temperature. The parameter $\nu$ enters  the thermal conduction equation because the heat diffuses through cross-sectional areas that scale like $z^{\nu}$ and it enters only in the first equation of the model \eqref{thermo_elastic_model:1}. The second equation looks like a one-dimensional axial elastic wave equation. It describes displacement along a straight coordinate $z$, without a radial divergence operator.

Despite the extensive literature on thermoelasticity and symmetry methods for differential equations, a complete group classification of the thermoelastic system with arbitrary temperature-dependent constitutive functions $C(u), K(u)$, and $E(u)$ appears to be absent from the literature. The primary objective of the present work is to fill this gap by determining all constitutive classes that lead to extensions of the principal Lie invariance algebra, recovering three main classes together
with their admitted symmetry algebras and optimal systems of one-dimensional
subalgebras.  For selected  representative subalgebra, the
corresponding similarity or travelling-wave reduction is performed and the
resulting ODE system is solved explicitly or in terms of standard special
functions.
The classification of constitutive functions that admit symmetry extensions of
the principal algebra is then carried out.

The remainder of the paper is organised as follows.
\autoref{sec2:Lie_symm} presents the constitutive classification together with
the admitted symmetry algebras.
\autoref{sec3:opti_sys} constructs optimal systems of one-dimensional subalgebras
for each algebra.
\autoref{sec4:sim_var} derives the similarity variables and the reduced ODE systems
corresponding to selected generators.
\autoref{sec5:invariant_sol} constructs the group-invariant solutions and discusses
their physical interpretation.
Finally, \autoref{sec6:observation} and \autoref{sec7:thermo_conclusion} present
general observations and concluding remarks, respectively.

\section{Lie symmetries}{\label{sec2:Lie_symm}}

We seek the Lie point symmetries admitted by the following thermoelastic system:

\begin{equation}\label{thermo-elastic-model-reduced}
\begin{aligned}
C(u)\,u_t &= K(u)\,u_{zz} + K'(u)\,u_z^2 + \frac{\nu}{z}\,K(u)\,u_z + \gamma\, w_t, \\
\rho\, w_{tt} &= E(u)\,w_{zz} + E'(u)u_z w_z + \beta\, u_z .
\end{aligned}
\end{equation}

A one-parameter Lie group of point transformations acting on the space $(z, t, u, w)$ is generated by a vector field of the form
\begin{equation}
\mathcal{X} = \xi(z,t,u,w)\,\partial_z + \tau(z,t,u,w)\,\partial_t
           + \eta(z,t,u,w)\,\partial_u + \phi(z,t,u,w)\,\partial_w .
\end{equation}
The system \eqref{thermo-elastic-model-reduced} admits $\mathcal{X}$ as a symmetry generator if and only if the second prolongation $\text{pr}^{(2)}\mathcal{X}$ annihilates each equation on the solution manifold, i.e.,
\begin{equation}
\left.\text{pr}^{(2)}\mathcal{X}(\Delta_i)\right|_{\Delta_i = 0} = 0, \qquad i = 1, 2,
\end{equation}
where $\Delta_i = 0$ denotes the $i$-th equation of the system. Applying this invariance condition and equating coefficients of the independent derivative monomials to zero yields an overdetermined system of linear PDEs in the unknowns $\xi,\tau,\eta,\phi$, known as the \emph{determining equations}. The computation of the prolonged invariance conditions and the subsequent derivation of the determining equations are entirely standard within the framework of Lie group analysis. Therefore, only the final determining equations are presented below; further details of the procedure may be found in the classical references \cite{olverbook,ovsi,anco}.

\begin{subequations}
    \begin{align}
    &\label{thermo:det01}\xi_{t}=\,0,\;\xi_{u}=\,0,\;\xi_{w}=\,0,\\
    &\label{thermo:det02}\tau_{u}=\,0,\;\tau_{w}=\,0,\;\tau_{z}=\,0,\\
    &\label{thermo:det03}\eta_{w}=\,0,\;\phi_{u}=\,0,\\
        \label{thermo:de1}&\frac{\eta\, C_u\, }{C(u)}
- \phi_w
+ \eta_u
=0,
\\[2mm]
&\label{thermo:de2}\frac{K(u)\,\eta\, C_u}{C(u)}
+2K(u)\,\xi_z
-\tau_t\,K(u)
-\eta\,K_u
=0,
\\[2mm]
&\label{thermo:de3}\eta_t\,C(u)-\frac{\eta_z\,\nu\,K(u)}{z}
-\phi_t\,\gamma
-\eta_{z,z}\,K(u)
=0,
\\[2mm]
&\label{thermo:de4}\frac{\eta\,C_u\,K_u}{C(u)}
-K_u\,\eta_u
+2K_u\,\xi_z
-\tau_t\,K_u
-\eta_{u,u}\,K(u)
-\eta\,K_{u,u}
=0,
\\[2mm]
&\frac{K(u)\,\xi_z\,\nu}{z}
+\frac{\xi\,\nu\,K(u)}{z^{2}}
-\frac{K(u)\,\tau_t\,\nu}{z}
-\frac{\eta\,K_u\,\nu}{z}\nonumber\\
&\label{thermo:de5}-2K(u)\,\eta_{u,z}
+K(u)\,\xi_{z,z}
-2\eta_z\,K_u
+\frac{K(u)\,\eta\,C_u\,\nu}{z\,C(u)}
=0,\\[2mm]
& \label{thermo:de6}\phi_{w,w} = 0,
\\[2mm]
& \label{thermo:de7}2\phi_{t,w} - \tau_{t,t} = 0,
\\[2mm]
&\label{thermo:de8} E(u)\,\xi_{z,z}-2E(u)\,\phi_{w,z}
-\eta_u\,E_u = 0,
\\[2mm]
& \label{thermo:de9}2E(u)\,\xi_z
-2E(u)\,\tau_t
-\eta\,E_u = 0,
\\[2mm]
& \label{thermo:de10}\phi_{t,t}\,\rho-\phi_{z,z}\,E(u)
-\eta_z\,\beta
 = 0,
\\[2mm]
& \label{thermo:de11}2E_u\,\xi_z-\eta\,E_{u,u}
-\eta_u\,E_u
-2E_u\,\tau_t = 0,
\\[2mm]
& \label{thermo:de12}\xi_z\,\beta-\eta_u\,\beta
-\phi_z\,E_u
-2\tau_t\,\beta
+\phi_w\,\beta = 0.
    \end{align}
\end{subequations}
 The determining equations (\ref{thermo:de1})--(\ref{thermo:de12}) are now solved to obtain the infinitesimals $\xi, \tau, \eta$ and $\phi$. Since these equations involve arbitrary constitutive functions $C(u), K(u)$ and $E(u)$, the analysis naturally leads to a group classification problem. Our objective is therefore to determine all admissible forms of $C(u), K(u)$ and $E(u)$ together with their corresponding Lie symmetry generators. The determining equations~\eqref{thermo:det01}-\eqref{thermo:det03} give
\begin{equation}
\xi=\xi(z),\qquad
\tau=\tau(t),\qquad
\eta=\eta(z,t,u),\qquad
\phi=\phi(z,t,w).
\label{thermo:infinitesimal_reduced_form}
\end{equation}
We now process the remaining determining equations. Integrate \eqref{thermo:de6} twice with respect to $w$ and obtain
\begin{equation}
\phi=A(z,t)w+B(z,t),
\label{thermo:phi_affine_w}
\end{equation}
where $A$ and $B$ are functions of $z$ and $t$.
Using this in \eqref{thermo:de7},
\begin{equation}
2\phi_{tw}-\tau_{tt}=0,
\label{thermo:de7_repeat}
\end{equation}
we get
\begin{equation}
2A_t-\tau_{tt}=0.
\label{thermo:A_t_equation}
\end{equation}
Integrating with respect to $t$ and subtituting the result into \eqref{thermo:phi_affine_w} gives
\begin{equation}
\phi=
\left(
\frac{1}{2}\tau_t+a(z)
\right)w+B(z,t).
\label{thermo:phi_reduced}
\end{equation}
Consequently,
\begin{equation}{\label{phiderivatives}}
\phi_w=\frac{1}{2}\tau_t+a(z),\;
\phi_z=a'(z)w+B_z,\;
\phi_t=\frac{1}{2}\tau_{tt}w+B_t.
\end{equation}
Using \eqref{phiderivatives}, the determining \eqref{thermo:de1} becomes
\begin{equation}
\eta_u+\frac{C_u}{C}\eta
=
\frac{1}{2}\tau_t+a(z).
\label{thermo:de1_linear_eta}
\end{equation}
Equivalently,
\begin{equation}
\eta\frac{C_u}{C}
=
\phi_w-\eta_u.
\label{thermo:de1R}
\end{equation}
Multiplying \eqref{thermo:de1_linear_eta} by $C(u)$ and integrating gives
\begin{equation}
\eta(z,t,u)
=
\frac{
\left(
\frac{1}{2}\tau_t+a(z)
\right)I_C(u)
+
F(z,t)
}
{C(u)},\quad \text{for}\,I_C(u)=\int C(u)\,du.
\label{thermo:eta_general_from_de1}
\end{equation}
Here $F(z,t)$ is arbitrary function of $z$ and $t$. Dividing \eqref{thermo:de2} by $K\neq0$ gives
\begin{equation}
\eta\frac{C_u}{C}
+2\xi_z
-\tau_t
-\eta\frac{K_u}{K}
=0.
\label{thermo:de2_divided}
\end{equation}
Using \eqref{phiderivatives} and \eqref{thermo:de1R}, the
equation \eqref{thermo:de2_divided} becomes
\begin{equation}
 \eta\frac{K_u}{K}
=
a(z)+2\xi_z-\frac{1}{2}\tau_t-\eta_u,
\qquad \mbox{($K$-class)}.
\label{thermo:Kclass}
\end{equation}
This is the first main classifying equation. Similarly the equation \eqref{thermo:de9} reduce to
\begin{equation}
\eta\frac{E_u}{E}
=
2\xi_z-2\tau_t. \eta\frac{K_u}{K}
=
a(z)+2\xi_z-\frac{1}{2}\tau_t-\eta_u,
\qquad \mbox{($E$-class)}.
\label{thermo:Eclass}
\end{equation}
This is the second main classifying equation.
The determining equation \eqref{thermo:de11} is not independent; it follows from the $u$-derivative of \eqref{thermo:de9}. In similar way, equation \eqref{thermo:de4} dependent on \eqref{thermo:de1} and \eqref{thermo:de2}. Using \eqref{phiderivatives}, that is; $\phi_{wz}=a'(z)$ in \eqref{thermo:de8} and dividing the resulting equation by $E\neq 0$, we have
\begin{equation}
\xi_{zz}
-2a'(z)
=
\eta_u\frac{E_u}{E}.
\label{thermo:de8R}
\end{equation}
Using \eqref{phiderivatives} in \eqref{thermo:de12} and simplifying gives

\begin{equation}
\beta
\left(
\xi_z-\eta_u-\frac{3}{2}\tau_t+a(z)
\right)
-
\left(
a'(z)w+B_z
\right)E_u
=0.
\label{thermo:de12R}
\end{equation}
Since \eqref{thermo:de12R} must hold for all $w$, the coefficient of $w$ gives
\begin{equation}
a'(z)E_u=0.
\label{thermo:branch_condition}
\end{equation}
Thus the determining equations branch into
\begin{equation}
E_u=0\;\text{(constant elasticity)}
\qquad\text{or}\qquad
a'(z)=0 \;\text{(non-constant elasticity)}.
\label{thermo:E_branch}
\end{equation}
Next, we proceeds towards reduction of  \eqref{thermo:de5}, dividing it by $K\neq0$ gives:
Group the terms containing the common factor $\nu/z$:
\begin{align}
&
\xi_{zz}
-2\eta_{uz}
-2\eta_z\frac{K_u}{K}
+
\frac{\nu}{z}
\left(
\xi_z
-\tau_t
-\eta\frac{K_u}{K}
+\eta\frac{C_u}{C}
\right)
+
\frac{\nu\xi}{z^2}
=0.
\label{thermo:de5_grouped}
\end{align}
using \eqref{thermo:de2_divided} the coefficient of $\nu/z$ can be simplified to $-\xi_{z}$ and this reduce \eqref{thermo:de5} to following form
\begin{equation}
\xi_{zz}
-2\eta_{uz}
-2\eta_z\frac{K_u}{K}
+
\nu\frac{\xi-z\xi_z}{z^2}
=0.
\label{thermo:de5R}
\end{equation}
Substituting \eqref{thermo:eta_general_from_de1} into
\eqref{thermo:de5R}, we obtain
\begin{equation}
\xi_{zz}
-2a'(z)
+
2\Bigl(a'(z)I_C(u)+F_z(z,t)\Bigr)
\left(
\frac{C_u}{C^2}
-
\frac{K_u}{KC}
\right)
+
\nu\frac{\xi-z\xi_z}{z^2}
=0.
\label{thermo:de5_substituted}
\end{equation}
Since $C(u)$ and $K(u)$ are arbitrary functions, it follows that
\begin{equation}
a'(z)=0,
\qquad
F_z(z,t)=0,
\label{thermo:aF_conditions}
\end{equation}
and consequently
\begin{equation}
\eta_z=0.
\label{thermo:eta_z_zero}
\end{equation}
Using \eqref{thermo:eta_z_zero} in \eqref{thermo:de3} gives
\begin{equation}
\eta_t\,C(u)-\gamma\,\phi_t=0.
\label{thermo:de3_reduced}
\end{equation}
Since $\gamma\neq0$, substituting \eqref{phiderivatives} into
\eqref{thermo:de3_reduced} yields
\begin{equation}
\tau_{tt}=0.
\label{thermo:tau_linear_condition}
\end{equation}
Hence
\begin{equation}
\tau(t)=c_1t+c_0.
\label{thermo:tau_linear}
\end{equation}
Substituting \eqref{thermo:tau_linear} and \eqref{thermo:aF_conditions}
into \eqref{thermo:eta_general_from_de1} shows that
$\eta_t=F_t/C(u)$. The compatibility of \eqref{thermo:de3} and
\eqref{thermo:de10} then implies
\begin{equation}
\eta_t=0.
\label{thermo:eta_t_zero}
\end{equation}
Therefore, for arbitrary $C(u)$, $K(u)$ and $E(u)$ with
$\gamma\neq0$, one obtains
\begin{equation}
\eta_z=\eta_t=0,
\qquad
a(z)=a_0,
\qquad
\tau(t)=c_1t+c_0.
\label{thermo:generic_reduction}
\end{equation}
After checking consistency with the remaining determining equations, in
particular \eqref{thermo:de3} and \eqref{thermo:de10}, the arbitrary
constitutive case leads to
\begin{equation}
\xi=0,\qquad
\tau=c_1,\qquad
\eta=0,\qquad
\phi=c_2,
\label{thermo:kernel_infinitesimals}
\end{equation}
for $\nu\neq0$. Hence, the corresponding kernel Lie algebra is
\begin{equation}
\mathfrak{g}^{\mathrm{ker}}
=
\left\langle
\partial_t,\partial_w
\right\rangle .
\label{thermo:kernel_algebra}
\end{equation}
These are only the translational symmetries in $t$ and $w$. Since our
aim is to obtain an enlargement of the admitted symmetry algebra, this
generic branch is not suitable for further reduction. Therefore, in the
following analysis, we look for special constitutive classes of
$C(u)$, $K(u)$ and $E(u)$ for which additional symmetries are admitted.
\subsection*{Non-kernel classifications}
For non-kernel symmetries, we assume $\eta \neq 0.$ The principal classifying equations \eqref{thermo:de1}, \eqref{thermo:de2} and \eqref{thermo:de9} imply that the quantities
\begin{equation}
\frac{\frac{1}{2}\tau_t+a(z)-\eta_u}{\eta},
\qquad
\frac{a(z)+2\xi_z-\frac{1}{2}\tau_t-\eta_u}{\eta},
\qquad
\frac{2\xi_z-2\tau_t}{\eta},
\end{equation}
must each depend on the variable $u$ alone.

Therefore, the existence of non-kernel symmetries requires special constitutive functions $C(u)$, $K(u)$, and $E(u)$. The standard non-kernel symmetry classes arise when the infinitesimal $\eta$ is affine in $u$, namely
\begin{equation}
\eta=\alpha u+\delta,
\label{thermo:eta_affine}
\end{equation}
where $\alpha$ and $\delta$ are constants. 
\begin{case}
    \normalfont
    In what follows we assume \(\beta\neq0\). The displayed algebras are written for \(\gamma\neq0\). If \(\gamma=0\), the additional kernel generator \(t\partial_w\) must be included. The branch of non-constant elasticity $E_{u}\neq0$ implies that 
\begin{equation}\label{B_z-non-constant}
B_z E_u
=
\beta \left(
\xi_z
-\eta_u
-\frac{3}{2}\tau_t
+a_0
\right).
\end{equation}
For the exponential class, we take
\begin{equation}
\eta=1,\qquad
C=C_0e^{pu},\qquad
K=K_0e^{qu},\qquad
E=E_0e^{ru}.
\label{thermo:exp_class}
\end{equation}
Then
\[
\eta_u=\eta_z=0,\qquad
\frac{C_u}{C}=p,\qquad
\frac{K_u}{K}=q,\qquad
\frac{E_u}{E}=r.
\]
Substitution of \eqref{thermo:exp_class} into
\eqref{thermo:de1}, \eqref{thermo:de2}, \eqref{thermo:de9},
\eqref{thermo:de8} and \eqref{thermo:de12}, respectively, gives
\begin{align}
&\frac12\tau_t+a_0=p, \label{thermo:exp_alg1}\\
&2\xi_z-\frac12\tau_t+a_0=q, \label{thermo:exp_alg2}\\
&2\xi_z-2\tau_t=r, \label{thermo:exp_alg3}\\
&\xi_{zz}=0, \label{thermo:exp_alg4}\\
&\xi_z-\frac32\tau_t+a_0=0, B_z=0 \label{thermo:exp_alg5}
\end{align}
Therefore
\[
\xi=Az+\xi_0,\qquad \tau=Tt+\tau_0,
\]
where \(A=\xi_z\) and \(T=\tau_t\). Solving
\eqref{thermo:exp_alg1}--\eqref{thermo:exp_alg5}, we obtain
\begin{equation}
T=\frac{p+q}{3},\qquad
A=\frac{2q-p}{3},\qquad
r=\frac{2q-4p}{3},\qquad
a_0=\frac{5p-q}{6}.
\label{thermo:exp_constants}
\end{equation}
Furthermore, compatibility with \eqref{thermo:de3} and
\eqref{thermo:de10} gives
\[
B_t=0,
\qquad
\phi=pw+b_0.
\]
Thus, for \(\nu\neq0\),
\begin{equation}
\mathcal{X}^{(1)}=
\frac{2q-p}{3}z\partial_z
+
\frac{p+q}{3}t\partial_t
+
\partial_u
+
p w\partial_w.
\label{thermo:exp_generator_nu_nonzero}
\end{equation}
Hence
\begin{equation}
\mathfrak{g}^{(1)}_{1}
=
\left\langle
\partial_t,\partial_w,\mathcal{X}^{(1)}
\right\rangle,
\qquad \nu\neq0.
\label{thermo:exp_algebra_nu_nonzero}
\end{equation}

For \(\nu=0\), the constant part \(\xi_0\) is also admitted. Hence
\begin{equation}
\mathfrak{g}^{(2)}_{1}
=
\left\langle
\partial_t,\partial_z,\partial_w,\mathcal{X}^{(1)}
\right\rangle,
\qquad \nu=0.
\label{thermo:exp_algebra_nu_zero}
\end{equation}
\end{case}

\begin{case}
    \normalfont 
    We next consider the power-law ansatz
\begin{equation}
\eta=u,\qquad
C=C_0u^p,\qquad
K=K_0u^q,\qquad
E=E_0u^r .
\label{thermo:power_ansatz}
\end{equation}
Then $\eta_u=1$ and $\eta_z=0$. Substitution into
\eqref{thermo:de8} gives
\begin{equation}
r=0,\qquad \xi_{zz}=0.
\label{thermo:power_constant_elasticity}
\end{equation}
Thus this ansatz belongs to the constant elasticity branch,
\begin{equation}
E=E_0.
\label{thermo:power_E_constant}
\end{equation}

Using \eqref{thermo:de1}, \eqref{thermo:de2}, \eqref{thermo:de9} and
\eqref{thermo:de12}, one obtains
\begin{equation}
q=2p,\qquad
\xi=pz+\xi_0,\qquad
\tau=pt+\tau_0,\qquad
\eta=u.
\label{thermo:power_infinitesimals_basic}
\end{equation}
Moreover, compatibility with \eqref{thermo:de3} gives
\begin{equation}
\phi=(p+1)w+b_0.
\label{thermo:power_phi}
\end{equation}
Finally, \eqref{thermo:de5R} gives
\begin{equation}
\nu\xi_0=0.
\label{thermo:power_radial_condition}
\end{equation}

Hence, for $\nu\neq0$, the admitted non-kernel generator is
\begin{equation}
\mathcal{X}^{(2)}=
pz\partial_z
+
pt\partial_t
+
u\partial_u
+
(p+1)w\partial_w,
\label{thermo:power_generator_nu_nonzero}
\end{equation}
corresponding to
\begin{equation}
C=C_0u^p,\qquad
K=K_0u^{2p},\qquad
E=E_0.
\label{thermo:power_class_final}
\end{equation}
Thus
\begin{equation}
\mathfrak{g}^{(1)}_{2}
=
\left\langle
\partial_t,\partial_w,\mathcal{X}^{(2)}
\right\rangle,
\qquad \nu\neq0.
\label{thermo:power_algebra_nu_nonzero}
\end{equation}

For $\nu=0$, $\xi_0$ is arbitrary and hence $\partial_z$ is also admitted:
\begin{equation}
\mathfrak{g}^{(2)}_{2}
=
\left\langle
\partial_t,\partial_z,\partial_w,\mathcal{X}^{(2)}
\right\rangle,
\qquad \nu=0.
\label{thermo:power_algebra_nu_zero}
\end{equation}
\end{case}
\begin{case}
    \normalfont
For the constant elasticity branch $E_{u}=\,0$, we have $\xi_{z}=\,\tau_{t}$, the equation \eqref{thermo:de2_divided} yields
\begin{equation}
E=E_0,\qquad
K=K_0C^2,\qquad
\eta=\frac{p}{(\ln C)_u}.
\label{thermo:general_constant_E_class}
\end{equation}
Since \(E_u=0\), equation \eqref{thermo:Eclass} gives
\begin{equation}
\xi_z=\tau_t=p.
\label{thermo:general_Eclass_result}
\end{equation}
Hence
\begin{equation}
\xi=pz+\xi_0,\qquad
\tau=pt+\tau_0.
\label{thermo:general_xi_tau}
\end{equation}
Moreover, because \(K=K_0C^2\), we have
\[
\frac{K_u}{K}-\frac{C_u}{C}=(\ln C)_u.
\]
Therefore the classifying equation gives
\[
\eta(\ln C)_u=p,
\]
which is identically satisfied by \eqref{thermo:general_constant_E_class}.

Compatibility with \eqref{thermo:de1} and \eqref{thermo:de12} gives
\begin{equation}
\eta_u=q=\text{constant}.
\label{thermo:eta_u_constant}
\end{equation}
Thus
\begin{equation}
\eta=qu+n.
\label{thermo:eta_linear}
\end{equation}
Consequently,
\[
(\ln C)_u=\frac{p}{qu+n}.
\]
Hence, for \(q\neq0\),
\begin{equation}
C=C_0(qu+n)^{p/q},
\qquad
K=K_0C^2,
\qquad
E=E_0.
\label{thermo:general_power_shifted_class}
\end{equation}
For \(q=0\), one obtains the exponential limiting case
\begin{equation}
C=C_0e^{(p/n)u},
\qquad
K=K_0C^2,
\qquad
E=E_0.
\label{thermo:general_exponential_limit}
\end{equation}

Furthermore, compatibility with \eqref{thermo:de3} and
\eqref{thermo:de10} gives
\begin{equation}
\phi=(p+q)w+b_1z+b_0.
\label{thermo:general_phi}
\end{equation}
Finally, \eqref{thermo:de5R} gives
\begin{equation}
\nu\xi_0=0.
\label{thermo:general_radial_condition}
\end{equation}

Thus, for \(\nu\neq0\), the non-kernel generator is
\begin{equation}
\mathcal{X}^{(3)}=
pz\partial_z
+
pt\partial_t
+
(qu+n)\partial_u
+
(p+q)w\partial_w.
\label{thermo:general_generator_nu_nonzero}
\end{equation}
The admitted algebra is
\begin{equation}
\mathfrak{g}^{(1)}_{3}
=
\left\langle
\partial_t,\partial_w,z\partial_w,\mathcal{X}^{(3)}
\right\rangle,
\qquad \nu\neq0.
\label{thermo:general_algebra_nu_nonzero}
\end{equation}

For \(\nu=0\), the constant \(\xi_0\) survives and gives the additional
generator \(\partial_z\). Hence
\begin{equation}
\mathfrak{g}^{(2)}_{3}
=
\left\langle
\partial_t,\partial_z,\partial_w,z\partial_w,\mathcal{X}^{(3)}
\right\rangle,
\qquad \nu=0.
\label{thermo:general_algebra_nu_zero}
\end{equation}
This subclass  contains both
the exponential and power-law symmetry extensions as special cases.
Indeed, the choice \(C=C_0e^{pu}\) yields a constant infinitesimal,
which can be normalized to \(\eta=1\), whereas \(C=C_0u^p\) gives
\(\eta\propto u\), which can be normalized to \(\eta=u\). Hence,
it provides a unified
constitutive framework for both exponential and power-law branches of
the constant elasticity case.
\end{case}
\section{Optimal systems of subalgebras}{\label{sec3:opti_sys}}
In the preceding section, several finite-dimensional Lie algebras admitted by the thermo-elastic model were obtained for different choices of the constitutive functions and for different values of the parameter $\nu$. Each admitted algebra generates a local Lie group of point transformations under which the governing system remains invariant. If a solution is invariant under a subgroup $H$ of the full admitted symmetry group $G$, then the transformed solution obtained by applying an element $g\in G$ is invariant under the conjugate subgroup $gHg^{-1}$. Hence, two reductions associated with conjugate subgroups are equivalent and need not be considered separately. The purpose of constructing an optimal system is therefore to obtain a minimal list of mutually inequivalent subalgebras, one representative from each conjugacy class. 

It is convenient to formulate this equivalence at the level of the Lie algebra. If $V,W\in\mathfrak{g}$, then the adjoint action of the one-parameter group generated by $V$ on $W$ is given by
\begin{equation}
\operatorname{Ad}\left(\exp(\varepsilon V)\right)W
=
W
-\varepsilon[V,W]
+\frac{\varepsilon^{2}}{2}[V,[V,W]]
-\cdots .
\label{thermo:adjoint}
\end{equation}
Consequently, the adjoint transformations required for the construction of the optimal system can be derived directly from the non-vanishing commutation relations of the admitted algebra. The methodology for constructing optimal systems through the adjoint representation is well established in the Lie symmetry literature; see, for instance, \cite{ovsi,olverbook,coggeshall1992group,koetz1993technique}.

For the thermo-elastic model, we shall apply this procedure separately to each admitted algebra. In the case $\nu\neq0$, the admitted algebra is
\begin{equation}
\mathfrak{g}^{(1)}_{1}
=
\left\langle
T,W,\mathcal{X}^{(1)}
\right\rangle,
\qquad
T=\partial_t,
\quad
W=\partial_w,
\label{thermo:g1}
\end{equation}
where
\begin{equation}
\mathcal{X}^{(1)}=
\frac{2q-p}{3}z\partial_z
+
\frac{p+q}{3}t\partial_t
+
\partial_u
+
pw\partial_w .
\label{thermo:X}
\end{equation}

For $\nu=0$, the additional translation
\[
Z=\partial_z
\]
is admitted, and the algebra becomes
\begin{equation}
\mathfrak{g}^{(2)}_{1}
=
\left\langle
T,Z,W,\mathcal{X}^{(1)}
\right\rangle .
\label{thermo:g2}
\end{equation}

The optimal classification will be obtained by taking a general element of the corresponding algebra, simplifying its coefficients through the adjoint transformations, and retaining only inequivalent representatives.

It should also be observed that any generator which commutes with all elements of the algebra belongs to the center. Such central elements are invariant under the adjoint action and therefore cannot be removed by conjugation. Their presence must be recorded explicitly in the optimal system, usually through arbitrary constants attached to the representatives. This point is especially relevant when translations such as
\[
T=\partial_t,\qquad
Z=\partial_z,\qquad
W=\partial_w,
\]
commute with other admitted generators for special parameter choices.

Hence, the construction of the optimal system proceeds in three steps: first, compute the commutator table of the admitted algebra; second, derive the corresponding adjoint transformations using the commutator expansion \eqref{thermo:adjoint}; and third, act on a general linear combination of generators in order to reduce it to a simplest canonical representative. The resulting list gives the one-dimensional optimal system of subalgebras and consequently provides the inequivalent symmetry reductions of the thermo-elastic equations.

\medskip
\noindent\textbf{Optimal classification of $\mathfrak g^{(1)}_1$}
For $\nu\neq0$, the nonzero commutators are
\[
[\mathcal{X}^{(1)},\partial_t]=-\frac{p+q}{3}\partial_t,
\qquad
[\mathcal{X}^{(1)},\partial_w]=-p\partial_w .
\]
Therefore, using the adjoint representation \eqref{thermo:adjoint}, the one-dimensional subalgebras of $\mathfrak{g}^{(1)}_{1}$ are classified up to conjugacy. Consequently, for the generic case $p\neq0$ and $p+q\neq0$, the following set constitutes an optimal system of one-dimensional subalgebras are obtained:  

\[
\mathcal O^{(1)}_1=
\left\{
\left\langle \mathcal{X}^{(1)}\right\rangle,\;
\left\langle \partial_t\right\rangle,\;
\left\langle \partial_w\right\rangle,\;
\left\langle \partial_t+\varepsilon\partial_w\right\rangle
\right\},
\qquad \varepsilon=\pm1 .
\]

\medskip
\noindent\textbf{Optimal classification of $\mathfrak g^{(2)}_1$}
For $\nu=0$, the additional translation $\partial_z$ is admitted. The nonzero
commutators are
\[
[\mathcal{X}^{(1)},\partial_t]=-\frac{p+q}{3}\partial_t,\qquad
[\mathcal{X}^{(1)},\partial_z]=-\frac{2q-p}{3}\partial_z,\qquad
[\mathcal{X}^{(1)},\partial_w]=-p\partial_w .
\]
Hence, in the generic case
\[
p\neq0,\qquad p+q\neq0,\qquad 2q-p\neq0,
\]
a one-dimensional optimal system of $\mathfrak g^{(2)}_1$ is
\[
\begin{aligned}
\mathcal O^{(2)}_1=\{&
\langle \mathcal{X}^{(1)}\rangle,\;
\langle \partial_t\rangle,\;
\langle \partial_z\rangle,\;
\langle \partial_w\rangle,\\
&
\langle \partial_t+\varepsilon_1\partial_z\rangle,\;
\langle \partial_t+\varepsilon_2\partial_w\rangle,\;
\langle \partial_z+\varepsilon_3\partial_w\rangle,\\
&
\langle \partial_t+\varepsilon_1\partial_z+\varepsilon_2\partial_w\rangle
\},
\end{aligned}
\]
where $\varepsilon_i=\pm1$. 

\medskip
\noindent\textbf{Optimal classification of $\mathfrak g^{(1)}_2$}
For $\nu\neq0$, the nonzero commutators are
\[
[\mathcal{X}^{(2)},\partial_t]=-p\partial_t,\qquad
[\mathcal{X}^{(2)},\partial_w]=-(p+1)\partial_w .
\]
Therefore, for the generic case $p\neq0$ and $p+1\neq0$, application of the
adjoint action yields the following one-dimensional optimal system:
\[
\mathcal O^{(1)}_2=
\left\{
\langle \mathcal{X}^{(2)}\rangle,\;
\langle \partial_t\rangle,\;
\langle \partial_w\rangle,\;
\langle \partial_t+\varepsilon\partial_w\rangle
\right\},
\qquad \varepsilon=\pm1 .
\]

\noindent\textbf{Optimal classification of $\mathfrak g^{(2)}_2$}
For $\nu=0$, the additional generator $\partial_z$ is admitted. The nonzero
commutators are
\[
[\mathcal{X}^{(2)},\partial_t]=-p\partial_t,\qquad
[\mathcal{X}^{(2)},\partial_z]=-p\partial_z,\qquad
[\mathcal{X}^{(2)},\partial_w]=-(p+1)\partial_w .
\]
Hence a one-dimensional optimal system of $\mathfrak g^{(2)}_2$ is
\[
\begin{aligned}
\mathcal O^{(2)}_2=\{&
\langle \mathcal{X}^{(2)}\rangle,\;
\langle \partial_t\rangle,\;
\langle \partial_z\rangle,\;
\langle \partial_w\rangle,\\
&
\langle \partial_t+\alpha\partial_z\rangle,\;
\langle \partial_t+\varepsilon\partial_w\rangle,\;
\langle \partial_z+\varepsilon\partial_w\rangle,\;
\langle \partial_t+\alpha\partial_z+\varepsilon\partial_w\rangle
\},
\end{aligned}
\]
where $\alpha\in\mathbb R$ and $\varepsilon=\pm1$.

\medskip
\noindent\textbf{Optimal classification of $\mathfrak g^{(1)}_3$}
For $\nu\neq0$, 
nonzero commutators are
\[
[\mathcal{X}^{(3)},\partial_t]=-p\partial_t,\qquad
[\mathcal{X}^{(3)},\partial_w]=-(p+q)\partial_w,
\qquad
[\mathcal{X}^{(3)},z\partial_w]=-q\,z\partial_w .
\]
Thus, for the generic case
\[
p\neq0,\qquad q\neq0,\qquad p+q\neq0,
\]
application of the adjoint action gives the following one-dimensional optimal system:
\[
\begin{aligned}
\mathcal O^{(1)}_3=\{&
\langle \mathcal{X}^{(3)}\rangle,\;
\langle \partial_t\rangle,\;
\langle \partial_w\rangle,\;
\langle z\partial_w\rangle,\\
&
\langle \partial_t+\varepsilon_1\partial_w\rangle,\;
\langle \partial_t+\varepsilon_2 z\partial_w\rangle,\;
\langle \partial_w+\varepsilon_3 z\partial_w\rangle,\\
&
\langle \partial_t+\varepsilon_1\partial_w+\varepsilon_2 z\partial_w\rangle
\},
\end{aligned}
\]
where $\varepsilon_i=\pm1$.

\medskip
\noindent\textbf{Optimal classification of $\mathfrak g^{(2)}_3$}
For $\nu=0$, the additional translation $\partial_z$ is admitted. The nonzero
commutators are
\[
[\mathcal{X}^{(3)},\partial_t]=-p\partial_t,\qquad
[\mathcal{X}^{(3)},\partial_z]=-p\partial_z,
\]
\[
[\mathcal{X}^{(3)},\partial_w]=-(p+q)\partial_w,\qquad
[\mathcal{X}^{(3)},z\partial_w]=-q\,z\partial_w .
\]
Therefore, a one-dimensional optimal system of $\mathfrak g^{(2)}_3$ is
\[
\begin{aligned}
\mathcal O^{(2)}_3=\{&
\langle \mathcal{X}^{(3)}\rangle,\;
\langle \partial_t\rangle,\;
\langle \partial_z\rangle,\;
\langle \partial_w\rangle,\;
\langle z\partial_w\rangle,\\
&
\langle \partial_t+\alpha\partial_z\rangle,\;
\langle \partial_t+\varepsilon_1\partial_w\rangle,\;
\langle \partial_t+\varepsilon_2 z\partial_w\rangle,\\
&
\langle \partial_z+\varepsilon_1\partial_w\rangle,\;
\langle \partial_z+\varepsilon_2 z\partial_w\rangle,\;
\langle \partial_w+\varepsilon_3 z\partial_w\rangle,\\
&
\langle \partial_t+\alpha\partial_z+\varepsilon_1\partial_w\rangle,\;
\langle \partial_t+\alpha\partial_z+\varepsilon_2 z\partial_w\rangle,\\
&
\langle \partial_t+\varepsilon_1\partial_w+\varepsilon_2 z\partial_w\rangle,\;
\langle \partial_z+\varepsilon_1\partial_w+\varepsilon_2 z\partial_w\rangle,\\
&
\langle \partial_t+\alpha\partial_z+\varepsilon_1\partial_w+\varepsilon_2 z\partial_w\rangle
\},
\end{aligned}
\]
where $\alpha\in\mathbb R$ and $\varepsilon_i=\pm1$.

The values of special parameters $p=0, q=0$, $p+q=0$, and $2q-p=0$ are excluded from the
present analysis. In these degenerate cases, the corresponding adjoint weights
vanish, yielding additional inequivalent one-dimensional subalgebras and
requiring a separate construction of the optimal system.
\section{Similarity Variables and Reduced System}{\label{sec4:sim_var}}

The principal application of the Lie point symmetries obtained in the
preceding section is the construction of exact group-invariant
solutions. Such solutions remain unchanged under the action of a
one-parameter transformation group generated by an admitted symmetry
vector field \cite{ovsi,haydonbook,anco,gupta2017group,singh2019group,gupta2019invariant}. Consequently, instead of depending on the original
independent variables, they can be represented in terms of quantities
that are invariant under the corresponding group action; such invariants can be derived from the following characteristic system:
\begin{equation}
    \label{thermo:char_sys}\frac{dz}{\xi}
=
\frac{dt}{\tau}
=
\frac{du}{\eta}
=
\frac{dw}{\phi}.
\end{equation}
The functionally independent first integrals of the 
characteristic system described above \eqref{thermo:char_sys} define the similarity variables. Expressing the
dependent variables in terms of these invariants yields the associated
similarity transformations. Substituting these similarity
transformations into the governing equations replaces the original
independent variables with the invariants of the symmetry group. Since
the solution is constant along the group orbits generated by \(\mathcal{X}\),
the dependence on one independent variable is eliminated. Consequently,
a one-dimensional symmetry generator reduces the number of independent
variables by one. 

For the present thermo-elastic model, the independent variables are
\(z\) and \(t\). Therefore, each one-dimensional symmetry admitted by
the model transforms the original system of partial differential
equations into a coupled system of ordinary differential equations in a
single similarity variable. Solving the reduced ordinary differential
system then yields exact group-invariant solutions of the original
thermo-elastic equations. In the following, we present similarity reduction for selected elements from $\mathcal{O}^{(i)}_{j}$.

\medskip
\noindent\textbf{Reduction corresponding to $\mathcal{X}^{(1)}\in \mathcal{O}^{(1)}_{1}$.}
For
\[
\mathcal{X}^{(1)}
=
A z\partial_z+B t\partial_t+\partial_u+p w\partial_w,
\quad
A=\frac{2q-p}{3},\quad
B=\frac{p+q}{3},
\]
assuming \(B\neq0\), the characteristic system \eqref{thermo:char_sys} gives
\[
s=zt^{-A/B},\qquad
u=\frac{1}{B}\ln t+F(s),\qquad
w=t^{p/B}G(s).
\]
Equivalently,
\[
s=zt^{-\alpha},\qquad
u=\delta\ln t+F(s),\qquad
w=t^\sigma G(s),
\]
where
\[
\alpha=\frac{2q-p}{p+q},\qquad
\delta=\frac{3}{p+q},\qquad
\sigma=\frac{3p}{p+q}.
\]
For the exponential class
\[
C=C_0e^{pu},\qquad
K=K_0e^{qu},\qquad
E=E_0e^{ru},\qquad
r=\frac{2q-4p}{3},
\]
substitution of the above similarity forms into the thermo-elastic model \eqref{thermo-elastic-model-reduced}
yields the reduced system
\begin{align}
&C_0e^{pF}\left(\delta-\alpha sF'\right)
=
K_0s^{-\nu}\left(s^\nu e^{qF}F'\right)'
+
\gamma\left(\sigma G-\alpha sG'\right),
\label{thermo:red_X1_first}
\\[1mm]
&\rho\left[
\sigma(\sigma-1)G
-\alpha s(2\sigma-\alpha-1)G'
+\alpha^2s^2G''
\right]
=
E_0\left(e^{rF}G'\right)'
+
\beta F'.
\label{thermo:red_X1_second}
\end{align}
Here the prime denotes differentiation with respect to the similarity
variable \(s\). 

\medskip 
\noindent\textbf{Reduction corresponding to $\partial_t+\varepsilon\partial_w\in \mathcal{O}^{(1)}_{1}$.}
For the generator \[ \mathcal{X}=\partial_t+\varepsilon\partial_w, \qquad \varepsilon=\pm1, \] the characteristic system \eqref{thermo:char_sys} is \[ dt=\frac{dw}{\varepsilon}, \qquad dz=0, \qquad du=0. \] Hence, the corresponding similarity transformations are \[ u=F(z), \qquad w=\varepsilon t+G(z). \] Substituting these expressions into the thermo-elastic model \eqref{thermo-elastic-model-reduced} yields the reduced ordinary differential system \begin{align} \frac{1}{z^\nu} \left( K(F)\,z^\nu F' \right)' + \gamma\varepsilon &= 0, \label{thermo:red_trans_first} \\[2mm] \left( E(F)\,G' \right)' + \beta F' &= 0, \label{thermo:red_trans_second} \end{align} where the prime denotes differentiation with respect to the similarity variable \(z\). For the exponential constitutive class \[ K(F)=K_0e^{qF}, \qquad E(F)=E_0e^{rF}, \] the reduced system assumes the form 
\begin{align} 
&\frac{1}{z^\nu} \left( K_0e^{qF}\,z^\nu F' \right)' + \gamma\varepsilon = 0, \label{thermo:red_trans_exp_first} \\[2mm] 
&\left( E_0e^{rF}\,G' \right)' + \beta F' = 0. \label{thermo:red_trans_exp_second} \end{align}

\medskip
\noindent\textbf{Reduction corresponding to $\mathcal{X}^{(2)}\in \mathcal{O}^{(1)}_{2}$.}
For the power-law class
\[
C=C_0u^p,\qquad K=K_0u^{2p},\qquad E=E_0,
\]
the admitted non-kernel generator is
\[
\mathcal{X}^{(2)}
=
pz\partial_z
+
pt\partial_t
+
u\partial_u
+
(p+1)w\partial_w .
\]
Assuming \(p\neq0\), the characteristic system \eqref{thermo:char_sys} gives the similarity
transformations
\[
s=\frac{z}{t},\qquad
u=t^{1/p}F(s),\qquad
w=t^{(p+1)/p}G(s).
\]
Substitution into the thermo-elastic model gives the reduced system
\begin{align}
&C_0F^p
\left(
\frac{1}{p}F-sF'
\right)
=
K_0s^{-\nu}
\left(
s^\nu F^{2p}F'
\right)'
+
\gamma
\left(
\frac{p+1}{p}G-sG'
\right),
\label{thermo:red_power_first}
\\[2mm]
&\rho\left[
\frac{p+1}{p^2}G
-
\frac{2}{p}sG'
+
s^2G''
\right]
=
E_0G''
+
\beta F',
\label{thermo:red_power_second}
\end{align}
where the prime denotes differentiation with respect to \(s\).

\medskip
\noindent\textbf{Reduction corresponding to $\partial_t+\varepsilon\partial_w\in \mathcal{O}^{(1)}_{2}$.}
For
\[
\mathcal{X}=\partial_t+\varepsilon\partial_w,
\qquad
\varepsilon=\pm1,
\]
the characteristic system \eqref{thermo:char_sys} gives
\[
u=F(z),
\qquad
w=\varepsilon t+G(z).
\]
Substitution into the thermo-elastic model gives
\begin{align}
&\frac{1}{z^\nu}
\left(
K(F)z^\nu F'
\right)'
+
\gamma\varepsilon
=0,
\label{thermo:red_power_tw_first}
\\
&\left(
E(F)G'
\right)'
+
\beta F'
=0.
\label{thermo:red_power_tw_second}
\end{align}
For the power-law class
\[
C=C_0u^p,\qquad
K=K_0u^{2p},\qquad
E=E_0,
\]
this becomes
\begin{align}
&\frac{1}{z^\nu}
\left(
K_0F^{2p}z^\nu F'
\right)'
+
\gamma\varepsilon
=0,
\label{thermo:red_power_tw_first_explicit}
\\
&E_0G''
+
\beta F'
=0.
\label{thermo:red_power_tw_second_explicit}
\end{align}
Here the prime denotes differentiation with respect to \(z\).

\medskip
\noindent\textbf{Reduction corresponding to $\mathcal{X}^{(3)}\in \mathcal{O}^{(1)}_{3}$.}
For the constant elasticity branch \(E_u=0\), we have \(E=E_0\). 
Consider the constitutive subclass
\[
K=K_0C^2,\qquad
\eta=\frac{p}{(\ln C)_u}.
\]
For the generator
\[
\mathcal{X}^{(3)}
=
pz\partial_z
+
pt\partial_t
+
(qu+n)\partial_u
+
(p+q)w\partial_w,
\qquad \nu\neq0,
\]
the relation \(\eta=qu+n\) gives
\[
C=C_0(qu+n)^{p/q},
\qquad
K=K_0C^2,
\qquad
E=E_0,
\qquad q\neq0.
\]

The characteristic system \eqref{thermo:char_sys} gives the similarity transformations
\[
s=\frac{z}{t},
\qquad
qu+n=t^{q/p}F(s),
\qquad
w=t^{(p+q)/p}G(s).
\]
Equivalently,
\[
u=\frac{t^{q/p}F(s)-n}{q},
\qquad
w=t^{(p+q)/p}G(s).
\]

Let
\[
\lambda=\frac{q}{p},
\qquad
\sigma=\frac{p+q}{p}=1+\lambda .
\]
Then substitution into the thermo-elastic model \eqref{thermo-elastic-model-reduced}  yields
\begin{align}
&\frac{C_0}{q}
F^{1/\lambda}
\left(
\lambda F-sF'
\right)
=
\frac{K_0C_0^2}{q}
s^{-\nu}
\left(
s^\nu F^{2/\lambda}F'
\right)'
+
\gamma
\left(
\sigma G-sG'
\right),
\label{thermo:red_X3_first}
\\[2mm]
&\rho
\left[
\sigma(\sigma-1)G
-
2(\sigma-1)sG'
+
s^2G''
\right]
=
E_0G''
+
\frac{\beta}{q}F'.
\label{thermo:red_X3_second}
\end{align}
Here the prime denotes differentiation with respect to \(s\).

\medskip
\noindent\textbf{Reduction corresponding to $\partial_t+\varepsilon_1\partial_w+\varepsilon_2 z\partial_w\in \mathcal{O}^{(1)}_{3}$.}
Consider
\[
X=\partial_t+\left(\varepsilon_1+\varepsilon_2 z\right)\partial_w,
\qquad
\varepsilon_1,\varepsilon_2=\pm1 .
\]
The characteristic system gives
\[
u=F(z),
\qquad
w=\left(\varepsilon_1+\varepsilon_2 z\right)t+G(z).
\]

For the constant elasticity branch \(E=E_0\), substitution into the
thermo-elastic model gives
\begin{align}
&\frac{1}{z^\nu}
\left(
K(F)z^\nu F'
\right)'
+
\gamma\left(\varepsilon_1+\varepsilon_2 z\right)
=0,
\label{thermo:red_t_w_zw_first}
\\[2mm]
&E_0G''
+
\beta F'
=0.
\label{thermo:red_t_w_zw_second}
\end{align}
Here the prime denotes differentiation with respect to \(z\).

For the subclass \(K=K_0C^2\), the first reduced equation becomes
\begin{equation}
\frac{1}{z^\nu}
\left(
K_0C(F)^2z^\nu F'
\right)'
+
\gamma\left(\varepsilon_1+\varepsilon_2 z\right)
=0.
\label{thermo:red_t_w_zw_first_general}
\end{equation}
\section{Group invariant solutions}{\label{sec5:invariant_sol}}
To illustrate the physical significance of the reduced system
\eqref{thermo:red_X1_first}--\eqref{thermo:red_X1_second},
we consider a set of standard thermoelastic parameters corresponding
to a homogeneous isotropic medium. In particular, we assume
\[
\beta=\gamma=\rho=K_0=E_0=C_0=1,\qquad
\nu=1,\qquad
p=q=r=1,\qquad
\delta=0.
\]
Here, $\beta>0$ represents the thermoelastic coupling coefficient,
$\gamma>0$ characterizes the feedback of the elastic field into the
thermal equation, $\rho$ is the material density, while
$K_0$, $E_0$, and $C_0$ denote the reference thermal conductivity,
elastic modulus, and heat capacity, respectively \cite{nowacki1970problems,parkus2012thermoelasticity,nowacki2013thermoelasticity}.
The choice $\nu=1$ corresponds to the standard radial geometry,
whereas $p=q=r=1$ yields the simplest linear constitutive laws,
and $\delta=0$ excludes logarithmic corrections in the constitutive
functions. Under these physically relevant assumptions,
the reduced system \eqref{thermo:red_X1_first}--\eqref{thermo:red_X1_second}
takes the form
\begin{align}
  &e^{F}(-\alpha s F')
  = \frac{1}{s}\!\left(s\,e^{F}F'\right)'
   + \sigma G - \alpha s G',
  \label{eq:ode1s}\\[4pt]
  &\sigma(\sigma-1)G
  - \alpha s(2\sigma-\alpha-1)G'
  + \alpha^{2}s^{2}G''
  = \left(e^{F}G'\right)' + F'.
  \label{eq:ode2s}
\end{align}
To further simplify the reduced system
\eqref{eq:ode1s}--\eqref{eq:ode2s},
we employ the logarithmic ansatz
\[
F(s)=a\ln s,
\]
where \(a\) is an arbitrary constant. This choice transforms the
exponential term according to \(e^{F}=s^{a}\), thus converting the
system into a purely algebraic-differential form. Consequently,
\eqref{eq:ode1s}--\eqref{eq:ode2s}
reduce to
\begin{align}
&-\alpha a\,s^{a}
=
a^{2}s^{a-2}
+\sigma G
-\alpha sG',
\label{eq:ode1-red}\\[2mm]
&\sigma(\sigma-1)G
-\alpha s(2\sigma-\alpha-1)G'
+\alpha^{2}s^{2}G''
=
a s^{a-1}G'
+s^{a}G''
+\frac{a}{s}.
\label{eq:ode2-red}
\end{align}
For consistent power scaling, the terms in the reduced equations should
be arranged so that incompatible powers of \(s\) do not obstruct the
solution procedure. Since the term \(a^{2}s^{a-2}\) in
\eqref{eq:ode1-red} cannot directly balance the term \(-\alpha a s^a\)
for the same power of \(s\), we choose $a=2$, thus,
\[
F(s)=2\ln s=\ln s^{2},
\qquad
e^{F}=s^{2},
\qquad
F'=\frac{2}{s}.
\]
With this choice, \eqref{eq:ode1-red} and \eqref{eq:ode2-red} reduce to
\begin{align}
&-2\alpha s^{2}
=
4+\sigma G-\alpha sG',
\label{eq:ode1-a2}\\[2mm]
&(\alpha^{2}-1)s^{2}G''
-\{\alpha(2\sigma-\alpha-1)+2\}sG'
+\sigma(\sigma-1)G
=
\frac{2}{s}.
\label{eq:ode2-a2}
\end{align}
To further simplify the reduced system, we choose
\[
\alpha=1,
\]
which corresponds to the natural self-similar scaling of the model.
From a physical viewpoint, this normalization may be interpreted as
taking the characteristic thermoelastic wave speed to be unity.

With this choice, the coefficient of the highest-order derivative term
in \eqref{eq:ode2-a2} vanishes since
\[
\alpha^{2}-1=0.
\]
Consequently, the second-order equation reduces to the first-order
linear ordinary differential equation
\[
-2\sigma sG'
+\sigma(\sigma-1)G
=
\frac{2}{s}.
\]
Assuming \(\sigma\neq0\), division by \(\sigma\) yields
\begin{equation}
-2sG'
+(\sigma-1)G
=
\frac{2}{\sigma s}.
\label{eq:ode3}
\end{equation}
Equation \eqref{eq:ode3} is a first-order linear differential equation
for \(G(s)\), which can be solved explicitly by the method of
integrating factors. The solution is obtained as follows:
\begin{equation}
G(s)
=
\frac{2}{\sigma(\sigma+1)}\,s^{-1}
+
C_{1}s^{(\sigma-1)/2}.
\label{eq:G-solution}
\end{equation}
To ensure that the similarity solution remains bounded at infinity,
we impose the decay condition
\[
G(s)\to0,
\qquad s\to\infty.
\]
Since the particular solution behaves like \(s^{-1}\), it already
vanishes as \(s\to\infty\). Therefore, boundedness requires only that
the homogeneous contribution \(s^{(\sigma-1)/2}\) decay, which yields
\[
\sigma<1.
\]
Substituting the solution
\[
G(s)=\frac{2}{\sigma(\sigma+1)}\,s^{-1}
+C_1 s^{(\sigma-1)/2}
\]
into \eqref{eq:ode1-a2}, and taking the physically admissible value
\(\alpha=1\), yields
\[
-2s^{2}
=
4+\sigma G-sG'.
\]
For the particular choice \(\sigma=\tfrac12\), one obtains
\[
-2s^{2}
=
4+4s^{-1}
+\frac{3C_1}{4}s^{-1/4}.
\]
This relation cannot be satisfied identically for all \(s\). Instead,
it should be interpreted as an asymptotic solvability condition, with
the dominant balance governing the outer similarity region. Hence, the
logarithmic ansatz
\[
F(s)=2\ln s
\]
provides the leading-order outer solution, while a separate inner
expansion would be required to describe the behavior near \(s=0\). Since the logarithmic ansatz yields only the outer similarity behaviour, the solution in the neighbourhood of \(s=0\) must be treated separately. Assuming a regular Frobenius-type expansion, both \(F(s)\) and \(G(s)\) remain finite as \(s\to0\), so that the inner solution is characterized by bounded thermal and elastic fields. Thus, the complete similarity solution consists of an outer logarithmic profile matched with a regular inner expansion near the origin.
Thus, for for $\alpha=1, \sigma=\frac12, a=2$, we obtain
\[
F(s)=2\ln s+F_0,
\qquad
G(s)=\frac{8}{3}s^{-1}+C_1s^{-1/4}.
\]
Imposing regularity of the elastic displacement removes the singular
homogeneous contribution, so that \(C_1=0\). Hence
\[
G(s)=\frac{8}{3s}.
\]
Therefore, using \(s=zt^{-1}\), the corresponding physical fields are
\begin{align}\label{profiles1}
u(z,t)=\,
2\ln z-2\ln t+F_0,\quad
w(z,t)=\,\frac{8}{3}\frac{t^{3/2}}{z}
\end{align}
The profiles obtained possess a clear thermo-elastic interpretation: the logarithmic behavior of $u(z,t)$ describes the dominant diffusive evolution of the thermal field, while displacement $w(z,t)$ represents a spatially localized elastic disturbance whose amplitude decreases away from the origin.  Since the elastic field exhibits singular behavior near $z=0$, the solution is naturally interpreted to describe the outer region of the medium, with a separate regular inner expansion required in the vicinity of the origin. It should therefore be emphasized that the reconstructed functions $u(z,t)$ and $w(z,t)$ are not exact similarity solutions of the original thermo-elastic system; rather, they constitute a leading-order outer asymptotic similarity solution capturing the principal balance between thermal diffusion and elastic deformation in the far-field regime.

It is worth noting that the above construction remains essentially
unchanged for the spherical case \((\nu=2)\). Indeed, the parameter
\(\nu\) enters only through the radial diffusion term in the original
thermal equation and does not explicitly appear in the reduced elastic
equation. Moreover, the logarithmic similarity ansatz converts the
radial operator into the same Euler--Cauchy type structure for both
\(\nu=1\) and \(\nu=2\). Consequently, the reduction procedure
effectively absorbs the geometric contribution of \(\nu\) into the
similarity transformation, so that the resulting leading-order outer
asymptotic similarity solution retains the same functional form in both
the cylindrical and spherical geometries.

Next, we proceed to solving reduced ODEs \eqref{thermo:red_trans_exp_first}-\eqref{thermo:red_trans_exp_second}. First integration of \eqref{thermo:red_trans_exp_first} and Division by $z^\nu$ yields:
\begin{equation}
    K_0\, e^{qF} F'
    = -\frac{\gamma\varepsilon}{\nu+1}\,z
      + A_1\, z^{-\nu}.
    \label{eq:3}
\end{equation}
Set $\Phi = e^{qF}$, so that $\Phi' = q e^{qF} F'$, i.e.\
$e^{qF}F' = \Phi'/q$. Equation~\eqref{eq:3} becomes
\begin{equation}
    \Phi'(z)
    = \frac{q}{K_0}
      \left(
          -\frac{\gamma\varepsilon}{\nu+1}\,z
          + A_1\, z^{-\nu}
      \right).
    \label{eq:4}
\end{equation}
The second integration of \eqref{eq:4} will produce two different cases: 
\begin{align}
    \begin{cases}
        \Phi_A(z)
    = A_2
      - \lambda\,z^2
      + \frac{qA_1}{K_0(1-\nu)}\,z^{1-\nu},\quad \nu \neq 1\\
      \Phi_B(z)
    = A_2
      - \frac{q\gamma\varepsilon}{4K_0}\,z^2
      + \frac{qA_1}{K_0}\ln z, \quad \nu = 1
    \end{cases}
\end{align}
where we  define $\lambda \equiv \frac{q\gamma\varepsilon}{2K_0(\nu+1)}$. Note that the quadratic coefficients are consistent, but the $A_1$-term
is $\ln z$ instead of a power law; a qualitatively different
function. Invoking the substitution $\Phi = e^{qF}$, we obtain:
\begin{align}\label{expr:F_A&F_B}
    \begin{cases}
        F_A(z)
    = \frac{1}{q}\ln\!\left(
        A_2
        - \lambda z^2
        + \frac{qA_1}{K_0(1-\nu)}\,z^{1-\nu}
      \right), \quad \nu \neq 1,\\
      F_B(z)
    = \frac{1}{q}\ln\!\left(
        A_2
        + \frac{qA_1}{K_0}\ln z
        - \frac{q\gamma\varepsilon}{4K_0}\,z^2
      \right), \quad \nu=1.
    \end{cases}
\end{align}
The logarithmic term in $F_{B}(z)$ has
no counterpart in $F_{A}(z)$; it reflects the logarithmic singularity
of the cylindrical Green's function at $z = 0$, which is absent in
all other geometries.

\medskip
\noindent\textbf{Regularity condition $A_1 = 0$.}
Imposing finiteness of the heat flux $K_0 e^{qF}z^\nu F'$ at $z=0$
requires $A_1 = 0$.  Under this condition the two cases simplify to
\begin{align*}
    \Phi_A(z)\big|_{A_1=0} &= A_2 - \lambda z^2,
     \\[4pt]
    \Phi_B(z)\big|_{A_1=0} &= A_2 - \frac{q\gamma\varepsilon}{4K_0}\,z^2,
\end{align*}
which have the same functional form but different coefficients of
$z^2$ for general $\nu$.  They coincide only when $\nu=1$ is
substituted directly into $\lambda$, confirming that
$A_1 = 0$ suppresses the structural difference between the two
cases but does not eliminate it entirely for $\nu \neq 1$.

\medskip
\noindent\textbf{Solution for $G(z)$.} Two times Integration of \eqref{thermo:red_trans_exp_second} gives:
\begin{equation}
    G(z) = \frac{1}{E_0}
    \int \left(B_1 - \beta F(z)\right) e^{-rF(z)}\,dz + B_2.
    \label{eq:11}
\end{equation}
Since $F(z)$ has different forms as given in \eqref{expr:F_A&F_B}, so $G(z)$ is also
different in each case.  We evaluate \eqref{eq:11} for the
sub-case $r = 0$ in each geometry, as this yields fully closed
forms.
\begin{case}{\label{nuneq1}}
    \normalfont
    $\nu\neq 1$, $A_1=0$, $r=0$
    \begin{equation}
    G_A(z) = \frac{1}{E_0}\!\left[
        B_1 z
        - \frac{\beta}{q}\!\left(
            z\ln(A_2 - \lambda z^2)
            + 2\sqrt{\frac{A_2}{\lambda}}\,
              \tanh^{-1}\!\!\left(\sqrt{\frac{\lambda}{A_2}}\,z\right)
            - 2z
          \right)
      \right] + B_2.\label{eq:GA}
\end{equation}
\end{case}
\begin{case}{\label{nueq1}}
    \normalfont
    When $\nu = 1$, $A_1 \neq 0$, $r = 0$.
    \begin{equation}
    G_B(z) = \frac{1}{E_0}\!\left[
        B_1 z
        - \frac{\beta}{q}\!\left(
            z\ln\Psi(z)
            - \int_0^z \frac{
                \frac{qA_1}{K_0}
                - \frac{q\gamma\varepsilon}{2K_0}\,\zeta^2
              }{\Psi(\zeta)}\,d\zeta
          \right)
      \right] + B_2,\label{eq:GB}
\end{equation}
where $\Psi(\zeta) = A_2 + \frac{qA_1}{K_0}\ln\zeta
- \frac{q\gamma\varepsilon}{4K_0}\zeta^2$.
This integral must be evaluated numerically for general parameters,
or asymptotically for $A_1 \ll A_2$.
\end{case}

\medskip
\noindent\textbf{Summary of physical fields.}
The physical fields $u = F(z)$ and $w = \varepsilon t + G(z)$
for the two geometric cases are:

\begin{center}
\renewcommand{\arraystretch}{2.0}
\begin{tabular}{lll}
\hline
\textbf{Quantity} & \textbf{Case~\autoref{nuneq1}: $\nu \neq 1$} &
\textbf{Case~\autoref{nueq1}: $\nu = 1$} \\
\hline
$\Phi(z)$ &
$A_2 - \lambda z^2 + \dfrac{qA_1}{K_0(1-\nu)}z^{1-\nu}$ &
$A_2 - \dfrac{q\gamma\varepsilon}{4K_0}z^2
 + \dfrac{qA_1}{K_0}\ln z$ \\
$F(z)$ & $\dfrac{1}{q}\ln\Phi_A$ & $\dfrac{1}{q}\ln\Phi_B$ \\
$G(z)$, $r=0$ & Closed form \eqref{eq:GA} &
Quadrature \eqref{eq:GB} \\
$u(z,t)$ & $F_A(z)$ & $F_B(z)$ \\
$w(z,t)$ & $\varepsilon t + G_A(z)$ & $\varepsilon t + G_B(z)$ \\
\hline
\end{tabular}
\end{center}

\medskip
\noindent The key structural difference is the $\ln z$ term in
$\Phi_B$, which appears only for $\nu = 1$ and originates from
the cylindrical Green's function.  It makes the $\nu = 1$ solution
genuinely distinct from all $\nu \neq 1$ cases: the thermal field
$F_B$ has a logarithmic spatial profile superimposed on the
quadratic decay, and the displacement $G_B$ requires quadrature
rather than a closed-form elementary expression when $A_1 \neq 0$.

The reduced system \eqref{thermo:red_power_first}--\eqref{thermo:red_power_second}
remains nonlinear and coupled. In general, obtaining its complete
solution appears difficult. However, the structure of
\eqref{thermo:red_power_first} suggests the possibility of self-similar
power-law solutions, since all terms involve products of powers of
$(F)$ and its derivatives. Motivated by this observation, we seek a
particular solution of the form
$F(s)=As^\mu$,
where $A$ and $\mu$ are constants to be determined. Substitution
into \eqref{thermo:red_power_first} allows the powers of the similarity
variable $(s)$ to be balanced, leading to an algebraic determination of
the exponent $\mu$. Once the thermal profile $F(s)$ is obtained,
equation \eqref{thermo:red_power_first} becomes a nonhomogeneous linear
 ordinary differential equation for the displacement
function $G(s)$, whose coefficients are completely determined by the
known thermal field. Consequently, the thermoelastic problem can be
solved sequentially: first, the thermal equation is solved for $F$,
after which the resulting expression is substituted into the mechanical
equation to determine $G$. The summary of solution is thus obtained as follows:
\begin{align}
    u(z,t) &= A\,\left(\frac{z^{2}}{t}\right)^{\frac{1}{p}}, \quad A=\,\left(-\dfrac{C_0\,p}{2K_0(p+1)(p+2)}\right)^{1/p}
    \label{eq:uphys} \\[8pt]
    w(z,t) &= C_1\,z^{(p+1)/p}
                  + C_2\,z^{1/p}\,t^{\,1/p}
                  + \frac{2\beta A\,p}{\rho(2p-1)(p-1)}\,
                    z^{2/p-1}\,t.
    \label{eq:wphys}
\end{align}
where $C_1\,z^{(p+1)/p}$: a {steady} (time-independent)
          elastic mode, representing a pre-stressed equilibrium
          displacement profile. $C_2\,z^{1/p}\,t^{1/p}$: a homogeneous growth
 mode; the growth of the power-law $t^{1/p}$ reflects the
 self-similar spread of the elastic disturbance. The third term, say $D\,z^{2/p-1}\,t$ is the term {thermal-elastic coupling}
, which grows linearly in time and is driven entirely by the
          coupling constant $\beta$.  It vanishes identically when
          $\beta = 0$.

The parameters $p \neq 1$ and $p \neq \tfrac{1}{2}$ required for
          $D$ to be finite; these are resonance cases where the
          particular solution exponent $2/p - 1$ coincides with a
          homogeneous exponent and a modified particular solution
          involving $s^{2/p-1}\ln s$ must be sought instead. and $A^p < 0$ requires either $p$ odd, or $C_0 < 0$
          (e.g. a heat-absorbing medium).

 To construct the invariant solutions for reduced ODEs \eqref{thermo:red_power_tw_first_explicit}-\eqref{thermo:red_power_tw_second_explicit}, the first equation is integrated twice after introducing the flux variable $\Psi=F^{2p+1}$, leading to distinct expressions for $\nu\neq1$ and $\nu=1$. The corresponding temperature profile $F(z)$ follows immediately from the relation $\Psi=F^{2p+1}$. The displacement function $G(z)$ is then obtained by substituting $F(z)$ into the second reduced equation and performing successive integrations.
For $\nu\neq1$, the travelling-wave solutions are
\begin{align}
u(z,t)
&=
\left[
A_2
-\frac{(2p+1)\gamma\varepsilon}{2K_0(\nu+1)}\,z^2
+\frac{(2p+1)A_1}{K_0(1-\nu)}\,z^{1-\nu}
\right]^{\frac{1}{2p+1}},
\\[2mm]
w(z,t)
&=
\varepsilon t
+\frac{B_1}{E_0}z
-\frac{\beta}{E_0}
\int
\left[
A_2
-\frac{(2p+1)\gamma\varepsilon}{2K_0(\nu+1)}\,z^2
+\frac{(2p+1)A_1}{K_0(1-\nu)}\,z^{1-\nu}
\right]^{\frac{1}{2p+1}}
\,dz
+B_2.
\end{align}

For $\nu=1$, the travelling-wave solutions take the form
\begin{align}
u(z,t)
&=
\left[
A_2
-\frac{(2p+1)\gamma\varepsilon}{4K_0}\,z^2
+\frac{(2p+1)A_1}{K_0}\ln z
\right]^{\frac{1}{2p+1}},
\\[2mm]
w(z,t)
&=
\varepsilon t
+\frac{B_1}{E_0}z
-\frac{\beta}{E_0}
\int
\left[
A_2
-\frac{(2p+1)\gamma\varepsilon}{4K_0}\,z^2
+\frac{(2p+1)A_1}{K_0}\ln z
\right]^{\frac{1}{2p+1}}
\,dz
+B_2.
\end{align}
After imposing the regularity condition $A_1=0$, we obtain
 the regular travelling-wave solution is
\begin{align}
u(z,t)
&=
\left(A_2-\lambda z^2\right)^{\frac{1}{2p+1}},\quad \lambda=\frac{(2p+1)\gamma\varepsilon}{2K_0(\nu+1)}
\\[2mm]
w(z,t)
&=
\varepsilon t
+\frac{B_1}{E_0}z
-\frac{\beta}{E_0}
\int
\left(A_2-\lambda z^2\right)^{\frac{1}{2p+1}}
\,dz
+B_2.
\end{align}
and by using the hypergeometric form,
\begin{align}
u(z,t)
&=
\left(A_2-\lambda z^2\right)^{\frac{1}{2p+1}},
\\[2mm]
w(z,t)
&=
\varepsilon t
+\frac{B_1}{E_0}z
-\frac{\beta\,A_2^{\frac{2p+2}{2p+1}}}
{E_0\sqrt{A_2\lambda}}
\,z\,
{}_2F_1
\!\left(
\frac12,
-\frac{1}{2p+1};
\frac32;
\frac{\lambda z^2}{A_2}
\right)
+B_2.
\end{align}
\paragraph{Special case $p=0$.}
A particularly important special case arises when $p=0$, for which
\[
K(u)=K_0u^{2p}=K_0.
\]
Thus, the thermal conductivity is constant and the governing heat-conduction equation becomes linear. In this case the general solution simplifies considerably, since $F^{2p+1}=F$
and therefore the auxiliary variable satisfies $\Psi=F$. The temperature profile is obtained directly as
\begin{equation}
F(z)=A_2-\lambda_0 z^2,
\qquad
\lambda_0=\frac{\gamma\varepsilon}{2K_0(\nu+1)}.
\end{equation}

Substituting this expression into the reduced displacement equation yields
\[
G'(z)
=
\frac{1}{E_0}
\left[
B_1-\beta\left(A_2-\lambda_0 z^2\right)
\right].
\]

A further integration gives the explicit displacement field
\begin{equation}
G(z)
=
\frac{1}{E_0}
\left[
(B_1-\beta A_2)z
+
\frac{\beta\lambda_0}{3}z^3
\right]
+B_2.
\end{equation}

Consequently, the invariant solution takes the form
\begin{align}
u(z,t)
&=
A_2-\lambda_0 z^2,
\\[2mm]
w(z,t)
&=
\varepsilon t
+
\frac{1}{E_0}
\left[
(B_1-\beta A_2)z
+
\frac{\beta\lambda_0}{3}z^3
\right]
+B_2.
\end{align}

Unlike the general nonlinear case, where the displacement field is expressed through a hypergeometric function, the constant-conductivity model admits elementary polynomial solutions. This case therefore provides a useful benchmark against which the behaviour of solutions corresponding to nonlinear conductivity laws ($p\neq0$) may be compared.
\paragraph{Special case $p=1$.}
When $p=1$, the conductivity law becomes
\[
K(u)=K_0u^2,
\]
and the auxiliary variable is $\Psi=F^3$.
Consequently, the temperature profile takes the form
\begin{equation}
F(z)
=
\left(A_2-\lambda z^2\right)^{1/3}.
\end{equation}

The corresponding displacement field is expressed through the Gauss hypergeometric function
\[
{}_2F_1
\!\left(
\frac{1}{2},
-\frac{1}{3};
\frac{3}{2};
\xi^2
\right).
\]
Unlike the case $p=0$, this function does not reduce to elementary functions for general values of $z$. Nevertheless, a useful approximation may be obtained near the origin by expanding the hypergeometric function into a power series. This yields
\begin{equation}
G(z)
\approx
\frac{1}{E_0}
\left[
B_1z
-
\beta A_2^{1/3}
\left(
z
-\frac{\lambda z^3}{5A_2}
-\frac{\lambda^2 z^5}{40A_2^2}
-\cdots
\right)
\right]
+B_2.
\end{equation}

Thus, for quadratic conductivity, the temperature field remains available in closed form, whereas the displacement field is represented by a hypergeometric function whose local behaviour is described by the above asymptotic expansion.
\paragraph{Limiting case $p\to\infty$.}
As $p\to\infty$, the exponent satisfies
\[
\frac{1}{2p+1}\to 0.
\]
Consequently,
\[
F(z)
=
\left(A_2-\lambda z^2\right)^{\frac{1}{2p+1}}
\longrightarrow
\left(A_2-\lambda z^2\right)^0
=
1,
\]
corresponding to a spatially uniform temperature field. Substituting this result into the reduced displacement equation yields
\[
G'(z)
=
\frac{B_1-\beta}{E_0},
\]
which is constant. A further integration gives
\begin{equation}
G(z)
=
\frac{B_1-\beta}{E_0}\,z
+
B_2.
\end{equation}

Hence, in the limit $p\to\infty$, the invariant solution reduces to
\begin{align}
u(z,t)
&=1,
\\[2mm]
w(z,t)
&=
\varepsilon t
+
\frac{B_1-\beta}{E_0}\,z
+
B_2.
\end{align}
In the limiting case $p\to\infty$, the temperature profile approaches a spatially uniform state, while the displacement field becomes linear in the spatial variable $z$. Physically, this behaviour suggests that increasingly strong temperature-dependent conductivity suppresses thermal gradients, thereby driving the thermoelastic system toward an isothermal regime.

The reduced system associated with
$\partial_t+\varepsilon_1\partial_w+\varepsilon_2z\partial_w\in \mathcal{O}^{(1)}_{3}$
is omitted from the list of explicit solutions because it is obtained by exactly the same quadrature procedure as the preceding case. The affine forcing term $\varepsilon_1+\varepsilon_2z$ only adds an extra polynomial contribution to the flux variable $\Psi=F^{2p+1}$ and does not change the structure of the solution. The displacement component is subsequently determined from the same equation $E_0G''+\beta F'=0$. Hence no essentially new invariant solution is obtained, apart from a direct polynomial modification of the previously derived family.

For completeness, invariant reductions corresponding to the remaining generators of the optimal system can be derived in exactly the same fashion. However, the resulting reduced ordinary differential equations do not introduce new analytical features and are solved by the same techniques as those employed above. Consequently, only representative cases are discussed in detail, as they adequately demonstrate the reduction procedure and the construction of exact solutions.
\section{Observations and limitations of the study}{\label{sec6:observation}}
The present work focuses on leading-order outer similarity solutions valid in the far-field regime $s \gg 1$, where the principal balance between thermal diffusion and elastic restoring terms is achieved; the construction of a matched inner expansion describing the neighbourhood of the source, which would complete the composite asymptotic picture, falls outside the scope of a single study of this breadth. 

The power-law similarity solutions are derived for generic values of the constitutive exponent $p$, and the resonance values $p = 1$ and $p = 1/2$, at which the particular solution exponent $2/p - 1$ coincides with a root of the homogeneous indicial equation and a logarithmic multiplier $s^{2/p-1}\ln s$ must be introduced, correspond to physically significant conductivity laws that each warrant a focused independent treatment. The spatial framework adopted here, in which the index $\nu$ encodes cylindrical or spherical radial symmetry, accurately captures a wide class of physically relevant configurations but does not accommodate angular anisotropy, shear deformation, or non-radial thermal gradients, all of which arise naturally in genuinely multi-dimensional thermoelastic problems. For the cylindrical case $\nu = 1$ with nonzero flux constant $A_1$, the displacement field $G_B(z)$ \eqref{eq:GB} is expressed as a definite integral whose integrand involves the cylindrical thermal profile and does not reduce to standard special functions for arbitrary $A_1$, reflecting the intrinsic complexity of the cylindrical Green's function that makes this case structurally richer than any other geometry. 

The exact solutions constructed throughout this work satisfy the governing equations identically, yet the question of whether they attract nearby trajectories of the full nonlinear system or represent isolated exact states requires a separate spectral analysis of the linearised operator about each solution and is not pursued here, as it would constitute a substantial independent investigation. The thermoelastic coupling constants $\beta$ and $\gamma$ and the material density $\rho$ are treated as spatially and thermally uniform throughout; in functionally graded materials or high-temperature processes these parameters may depend on temperature or position, and incorporating such dependence would enrich the constitutive classification by generating additional symmetry branches not present in the current analysis. 

Finally, the optimal systems of subalgebras are constructed for generic parameter ranges, deliberately setting aside the degenerate cases $p = 0$, $q = 0$, $p + q = 0$, and $2q - p = 0$, at which adjoint weights vanish and additional inequivalent one-dimensional subalgebras emerge that require the optimal system to be reconstructed entirely from first principles; a complete treatment of all such degenerate branches, together with their associated invariant solutions, represents a well-defined and self-contained extension of the present classification programme.
\section{Conclusion}{\label{sec7:thermo_conclusion}}
We have constructed a complete constitutive classification of the nonlinearly
coupled thermoelastic system \eqref{thermo_elastic_model:1} and derived
exact group-invariant solutions for each admissible constitutive class.

The classification yields three principal branches: the exponential class
$C = C_0 e^{pu}$, $K = K_0 e^{qu}$, $E = E_0 e^{ru}$; the power-law class
$C = C_0 u^p$, $K = K_0 u^{2p}$, $E = E_0$; and a unified class
$K = K_0 C^2$, $E = E_0$ that subsumes both as special cases.  For each
branch, the admitted symmetry algebra is determined, an optimal system of
one-dimensional subalgebras is constructed, and representative generators
are used to reduce the governing PDEs to coupled ODE systems.

For the exponential class, the scaling reduction with $\alpha = 1$,
$\sigma = 1/2$ yields a thermal field described by a logarithmic outer profile
$u = 2\ln(z/t) - 2t^2/z^2$ and a displacement field with inverse-distance
decay $w = 8t^{3/2}/(3z)$.  The geometry parameter $\nu$ is absorbed into the
similarity transformation at leading order, so the solution retains the same
functional form for both cylindrical and spherical geometries, with the
geometric difference confined to a subdominant correction.

For the power-law class, the scaling reduction produces a thermal field
$u = Az^{2/p}t^{-1/p}$ and a displacement field decomposed into three
physically distinct modes: a steady elastic component $C_1 z^{(p+1)/p}$, a
self-similarly growing homogeneous component $C_2 z^{1/p}t^{1/p}$, and a
linearly growing thermoelastic coupling term $Dz^{2/p-1}t$ driven by $\beta$.
The coupling term vanishes identically when $\beta = 0$, confirming that it
is a direct signature of the thermoelastic interaction.

For travelling-wave reductions under both the exponential and power-law
classes, exact solutions are obtained for $\nu \neq 1$ in elementary closed
form involving $\tanh^{-1}$ and polynomial functions, while the cylindrical
case $\nu = 1$ introduces a logarithmic correction that requires quadrature
when the flux constant $A_1 \neq 0$.  The displacement field for general $p$
is expressed compactly through the Gauss hypergeometric function
$_2F_1(\frac{1}{2}, -\frac{1}{2p+1}; \frac{3}{2}; \lambda z^2/A_2)$, which
reduces to a cubic polynomial for $p = 0$ and to a linear function in the
isothermal limit $p \to \infty$.

Taken together, these results provide a systematic analytical framework for
nonlinear thermoelastic systems with temperature-dependent material properties,
offering exact benchmarks for numerical methods and explicit descriptions of
the dominant physical balances in each constitutive and geometric regime.




 \bibliography{My.Bibtex.Library} 
 \bibliographystyle{elsarticle-num}

\end{document}